\documentclass[letterpaper, 10 pt, conference]{ieeeconf}  

\IEEEoverridecommandlockouts                              

\usepackage{amsmath} 
\usepackage{amssymb}  
\usepackage{bm,upgreek}
\usepackage{dsfont}
\usepackage{mathtools}
\usepackage[yyyymmdd,hhmmss]{datetime}
\usepackage{ntheorem}
\usepackage{color}
\usepackage{epsfig} 
\usepackage{cite}
\usepackage{graphicx}
\usepackage{graphics}
\usepackage{epstopdf}
\mathtoolsset{showonlyrefs}
\newcommand{\Exp}[1]{\mathbb{E}[#1]}
\newcommand{\Prob}[1]{\mathbb{P}(#1)}          
 
\newtheorem{Assumption}{\textit{Assumption}}
\newtheorem{Lemma}{Lemma}
\theoremstyle{remark}
\newtheorem{Remark}{Remark}
\theoremstyle{Definition}
\newtheorem{Problem}{Problem}
\newtheorem{Theorem}{Theorem}

\title{\LARGE \bf
Team Optimal Control of Coupled Major-Minor Subsystems with Mean-Field Sharing 
}

\author{Jalal Arabneydi and Aditya Mahajan
\thanks{This work was supported by the
Natural Sciences and Engineering Research Council of Canada through Grant NSERC-RGPIN 402753-11.}
\thanks{Jalal Arabneydi and Aditya Mahajan are with Department of Electrical Engineering, 
        McGill University, 3480 University St., Montreal, QC, Canada.
        {\tt\small jalal.arabneydi@mail.mcgill.ca} and
        {\tt\small aditya.mahajan@mcgill.ca}}%
}

\begin{document}
\maketitle

\vspace*{-5.2cm}{\footnotesize{Proceedings of IEEE  Indian Control Conference, 2015.}}
\vspace*{4.45cm}
\thispagestyle{empty}
\pagestyle{empty}

\begin{abstract}
In this paper, we investigate team optimal control of coupled major-minor subsystems with mean-field sharing. In such a model, there is one major subsystem that directly influences the dynamics of $n$ homogeneous minor subsystems; however, the minor subsystems influence the dynamics of the major subsystem and each other only through their mean behaviour (indirectly). In this model, the major and the minor subsystems are arbitrarily coupled in the cost and have mean-field sharing information structure. We propose a two-step approach. In the first step, we describe a mean-field sharing model with multiple types of minor subsystems as a generalization of the mean-field sharing model of Arabneydi and Mahajan, CDC 2014. In the second step, we use the results obtained in the first step to construct a dynamic programming decomposition to identify optimal strategies for the major and the minor subsystems. We present an example with numerical results to illustrate the approach.
\end{abstract}

\section{Introduction}
 \subsection{Motivation}\label{Motivation}
 
Team theory---which is motivated by applications in networked control systems, economics, communication networks,  sensor networks, robotics, and transportation networks, etc. --- investigates  the optimal design of decentralized systems consisting of coupled subsystems controlled by multiple controllers that have access to different information but need to collaborate to accomplish a common task.

Research in team theory began with with the seminal work of Marschack and Radner \cite{MarschackRadner1972} in the context of organizational behaviour. Dynamic team problems have been investigated in systems and control following the seminal work of Witsenhansuen \cite{Witsenhausen1971separation, Witsenhausen1971information} and others  on dynamic team problems. The general dynamic team problem belongs in the NEXP complexity class \cite{Bernstein2002complexity}. As such, most of the research has focused on specific information structures that are amenable to analysis. We refer the reader to \cite{Yuksel2013stochastic, Mahajan2012Information} for  detailed surveys.

Another approach to analyzing large-scale decentralized systems is mean-field approximation which arose in statistical physics, specially in applications  consisting of a large number of interacting particles \cite{Balescue1975, Klimontovich1987s}. In such large-scale systems, a single particle has negligible effect on others and is influenced by other particles only through their mean effect. This property may be exploited to reduce many body problems into one body problems. This idea from statistical physics has been introduced in dynamical systems  by Huang, Caines,
and Mahlame \cite{Huang2003,HuangPeter2006, Nourian2013MM} and Lasry and Lions \cite{Lasry2007mean, LasryLions2011}. Similar ideas have also been used in large population games \cite{Weintraub2005oblivious}. We refer the reader to \cite{Gomes2013survey} for a survey.

In real world scenarios, there are many applications where one big (major) subsystem can directly influence many small (minor) subsystems while small subsystems can only influence the big subsystem by their average (collective) behaviour.  For example, consider a service-provider as the major subsystem and its users as the minor subsystems. The service-provider may directly influence each user by the cost of service. On the other hand, the average number of happy or unhappy users (collective behaviour) may affect the service-provider's reputation (and hence, its utility). In this paper, we study team optimal control of such systems  by generalizing our previous results on mean-field sharing \cite{Jalal2014MF}.

 The rest of this paper is organized as follows. In Section \ref{MF_models}, we present three models of mean-field sharing: \textit{basic mean-field sharing}, \textit{mean-field sharing with major-minor} subsystems, and \textit{mean-field sharing with multiple-types}. For brevity, in the sequel, we call these models the basic MF, MF-MM, and MF-T models, respectively. The reason for considering the MF-T model is that the MF-MM model may be viewed as a special case of the MF-T model; hence the main results of the MF-MM model can be obtained from the main results of the MF-T model. In Section \ref{Main_results}, we state the main results of the basic MF model in \cite{Jalal2014MF} and extend the solution methodology to the MF-T  and  MF-MM models.  In Section \ref{Example}, we present an example of the MF-MM model and in Section \ref{Conclusion}, we conclude this paper.

%

\subsection{Notation}

We use upper-case letters to denote random variables (e.g. $X$) and lower-case letters to denote their realizations (e.g. $x$). Short-hand notation $X_{a:b}$ is used for the vector $(X_a,X_{a+1},\ldots,X_{b})$ and bold letter is used to denote vectors e.g. $\mathbf Y=(Y^1,\ldots,Y^n)$ where $n$ is the size of vector $\mathbf Y$.  $\mathbb{P}(\cdot)$ is probability of an event, $\mathds{1}(\cdot)$ is indicator function of a set, and  $\mathbb{E}[\cdot]$ is expectation of a random variable. $\mathbb{N}$ and $\mathbb{R}^+$ refer to the set of natural numbers and   positive real numbers, respectively.

\section{Mean Field Models}\label{MF_models}

\subsection{Basic Mean-Field Sharing Model}\label{BMF_model}

We start with the description of the mean-field sharing model presented in \cite{Jalal2014MF}. We refer to this model as the basic mean-field (MF) model.

%

Consider a discrete time decentralized control system with $n \in \mathbb{N}$ homogeneous subsystems. Let $X^i_t \in \mathcal{X}$ and $U^i_t \in  \mathcal{U}$ denote the state and the control action of subsystem~$i$, $i \in \{1,\ldots,n\}$, at time $t, t \in \mathbb{N},$ respectively. $\mathcal{X}$ and $\mathcal{U}$ are  finite sets that do not depend on $i$. The \emph{mean-field} of the system  refers to the empirical distribution of the state of all subsystems. The mean-field  at time $t$, denoted by $Z_t$, is given by
\begin{equation}\label{BF_mean_field}
Z_t=\frac{1}{n} \sum_{i=1}^n \delta_{X^i_t}
\end{equation}
where $\delta_x$ denotes a Dirac measure on $\mathcal{X}$ with a point mass at $x$. Subsystem $i$ evolves as follows:
\begin{equation}\label{BMF_dynamics}
X^i_{t+1}=f_t(Z_t,X^i_t,U^i_t,W^i_t), \hspace{.2cm} i \in \{1,\ldots,n\}
\end{equation}
where $f_t$  denotes the dynamics at time $t$ (that do not depend on $i$) and $\{W^i_t\}_{t=1}^T$,  $i \in \{1,\ldots,n\}$, are independent processes. The \textit{primitive random variables} $(X^1_1,\ldots,X^n_1,\{W^1_t\}_{t=1}^T,\ldots,\{W^n_t\}_{t=1}^T)$ are mutually independent and defined on a common probability space. The probability distribution on $X^i_1$ and $W^i_t$ do not depend on $i$.

At time $t$, controller $i$ observes the history of the  mean-field $Z_{1:t}$ and its own local state $X^i_t$. Hence, the information set of controller $i$, $i \in \{1,\ldots,n\}$, at time $t$ is
\begin{equation}\label{Information_structure}
I^i_t=\{Z_{1:t},X^i_t\}.
\end{equation}
Controller $i$ chooses the control action $U^i_t \in \mathcal{U}$ according to
\begin{equation}\label{BMF_control_law}
U^i_t=g^i_t(Z_{1:t},X^i_t)
\end{equation}
where $g^i_t$ is the \emph{control law} of controller $i$.

To simplify the solution and ensure fairness between all controllers, we assume the following.
\begin{Assumption} \label{Homogeneous}
At any time $t$, the control laws at all controllers are identical i.e. $g^i_t=g^j_t$ for any $i,j \in \{1,\ldots,n\}$. Therefore, we drop the superscripts and denote the control law at every controller at time $t$ as $g_t$.
\end{Assumption}

As illustrated in \cite{Jalal2014MF}, in general,  Assumption \ref{Homogeneous} leads to a loss in performance. Nonetheless, it is a standard assumption in the literature on large-scale systems for reasons of simplicity, fairness, and robustness.

The system runs for a finite time horizon $T$, $T \in \mathbb{N}$. We refer to the collection of control laws $\mathbf{g}:=(g_1,\ldots,g_T)$ over time $T$  as the \textit{control strategy} of the system. Let $\mathcal{G}_{\mathit{BMF}}$ be a set of all such control strategies.

At time $t$, the system incurs a cost that depends on joint state $\mathbf X_t:=(X^1_t,\ldots,X^n_t)$ and joint action $\mathbf U_t:=(U^1_t,\ldots,U^n_t)$ as follows 
\begin{equation}\label{BMF_cost}
\ell_t(\mathbf X_t, \mathbf U_t).
\end{equation}
The performance of any strategy $\mathbf{g}$ is measured by the following expected total cost
\begin{equation}\label{BMF_total_cost}
 J(\mathbf{g})=\mathbb{E}^{\mathbf{g}}\left[\sum_{t=1}^T  \ell_t(\mathbf X_t, \mathbf U_t) \right]
\end{equation}
where the expectation is with respect to a joint measure induced on all system variables by the choice of $\mathbf{g}$. We are interested in the following optimization problem.
\begin{Problem}\label{BMF_problem}
For the basic MF model described above, identify a control strategy $\mathbf{g}^\ast$ to minimize the expected total cost $J(\mathbf{g})$ given by \eqref{BMF_total_cost}. 

\end{Problem}

In \cite{Jalal2014MF}, the basic MF model was analyzed using the common information approach of \cite{Nayyar2013CIA}. It was shown that the mean-field $Z_t$ is an information state for an appropriately defined coordinated system. Using this information state, an appropriate dynamic programming decomposition to identify an optimal control strategy was derived.

\subsection{Mean-Field Sharing with Major-Minor Subsystems (MF-MM) Model }\label{MMMF_model}

In this paper, we investigate a setup where, in addition to $n$ homogeneous subsystems (minor subsystems), there is one major subsystem that directly affects the evolution of all minor subsystems and the cost. We call this model, \textit{mean-field sharing with major and minor subsystems} (MF-MM for short). This model is motivated by the scenarios described in Section ~\ref{Motivation}. In mean-field games, a similar model was introduced by \cite{Huang2010lMM} and further analyzed in \cite{Nourian2013MM}.

 As before, let $X^i_t \in \mathcal{X}$ and $U^i_t \in \mathcal{U}$ denote the state and the control action of minor subsystem $i, i \in \{1,\ldots,n\}$, at time $t$, $t \in \mathbb{N}$. In addition, let $X^0_t \in \mathcal{X}^0$ and $U^0_t \in \mathcal{U}^0$ denote the state and the control action of the \textit{major subsystem}  at time $t$. The sets $\mathcal{X}^0$, $\mathcal{X}$, $\mathcal{U}^0$, and $\mathcal{U}$ are finite. Let $Z_t=\frac{1}{n} \sum_{i=1}^n \delta_{X^i_t}$ denote the mean-field of the \textit{minor subsystems}. The dynamics of the major and minor subsystems are given by 
\begin{equation}
\begin{aligned}
&X^0_{t+1}=f^0_t(Z_t,X^0_t,U^0_t,W^0_t),\\
&X^i_{t+1}=f_t(Z_t,X^0_t,X^i_t,U^i_t,W^i_t), \quad i \in \{1,\ldots,n\}
\end{aligned}
\end{equation}
where $\{W^i_t\}_{t=1}^T$,  $i \in \{1,\ldots,n\}$, are independent processes. 
As in the basic model, the primitive random variables $(X^0_1,X^1_1,\ldots,X^n_1,\{W^0_t\}_{t=1}^T,\{W^1_t\}_{t=1}^T,\ldots,\{W^n_t\}_{t=1}^T)$ are mutually independent. The dynamics $f_t$ and the probability distribution of $X^i_1$ and $W^i_t$ do not depend on $i$.

Note that, in addition to the mean-field, the dynamics of the minor subsystems are also influenced by the state of the major subsystem. The information sets of the major and  minor subsystems are
\begin{equation}\label{MF-MM_information_structure}
I^0_t=\{Z_{1:t},X^0_{1:t}\}, \quad 
I^i_t=\{Z_{1:t},X^0_{1:t},X^i_t\}.
\end{equation}
Thus, the history of the mean-field and the major subsystem's state are observed by both the major and minor controllers. In addition, the minor controller observes its local state as well.

As in the basic model, we assume that the control laws of all minor controllers are identical. Thus, the control actions are generated as follows
\begin{equation}\label{MF-MM_control_law}
U^0_t=g^0_t(Z_{1:t},X^0_{1:t}),\quad 
U^i_t=g_t(Z_{1:t},X^0_{1:t},X^i_t).
\end{equation}
The collection $(\mathbf{g}^0, \mathbf{g})$, where $\mathbf{g}^0=(g^0_1,\ldots,g^0_T)$ and $\mathbf{g}=(g_1,\ldots,g_T)$, is called the \textit{control strategy} of the system. Let $\mathcal{G}_{\mathit{MF-MM}}$ be a set of all such control strategies.

 At time $t$, system incurs a cost $ \ell_t(\mathbf{X}_t, X^0_t, \mathbf{U}_t, U^0_t)$ that depends on $\mathbf{X}_t, X^0_t, \mathbf{U}_t, U^0_t$. The performance of any strategy  $(\mathbf{g}^0, \mathbf{g})$ is quantified by the following expected total cost
\begin{equation}\label{MMMF_total_cost}
J(\bm g^0\hspace{-.1cm}, \bm g)\hspace{-.1cm}=\hspace{-.1cm}\mathbb{E}^{(\bm {g}^0\hspace{-.1cm},\bm g)} \hspace{-.15cm} \left[ \sum_{t=1}^T  \ell_t(\mathbf{X}_t, X^0_t, \mathbf{U}_t, U^0_t) \right]
\end{equation}
where the expectation is with respect to a joint measure induced on all system variables by the choice of $(\mathbf{g}^0, \mathbf{g})$. We are interested in the following optimization problem.

\begin{Problem}
For the MF-MM model, identify a control strategy $({\mathbf{g}^0}^\ast, \mathbf{g}^\ast)$ to minimize the expected total cost given by~\eqref{MMMF_total_cost}.
\end{Problem}

As mentioned earlier in this section, similar generalization of the basic MF model in the context of mean-field games (MFG) has been considered in \cite{Huang2010lMM,Nourian2013MM}\footnote{In MFG, the subsystems are only  allowed to be weakly coupled in the cost (through the mean-field $Z_t$ and state of the major subsystem $X^0_t$) while in this work subsystems are allowed to be arbitrarily coupled in the cost.}. Introducing a major agent gives rise to a non-trivial conceptual difficulty in MFG. To explain this difficulty, let us describe the solution approach used in MFG. Assume that all players follow a control strategy $\bm g=(g_1,\ldots,g_T)$, where $g_t: \mathcal{X} \mapsto \mathcal{U}$. In the case of infinite players, this strategy gives rise to a \textit{deterministic} trajectory of the mean-field $\{Z_t\}_{t=1}^T$. Given a deterministic mean-field, the problem of finding the best response at an agent is a centralized stochastic control problem. The corresponding optimal strategy is called the best response to $\{Z_t\}_{t=1}^T$. A strategy-mean field pair $(\bm g, \{Z_t\}_{t=1}^T)$ is called \textit{consistent} if $\{Z_t\}_{t=1}^T$ is generated by $\bm g$ and $\bm g$ is the best response to $\{Z_t\}_{t=1}^T$. 

In the major-minor setup described above, the mean-field generated by a strategy $\bm g=(g_1,\ldots,g_T)$, where $g_t: \mathcal{X}^0 \times \mathcal{X} \mapsto \mathcal{U}$, depends on the sample realization of $X^0_{1:T}$; hence, it is not deterministic. Thus, one can not directly use the solution approach described above.

In this paper, we generalize the results in \cite{Jalal2014MF} to the major-minor model; in particular we obtain a dynamic programming decomposition to determine an optimal strategy. In contrast to the conceptual difficulties that arise in the solution methodology of MFG, our solution is straightforward. 

We first describe a mean-field model with multiple types (MF-T for short, explained below) and show that an appropriate change of variables reduces the MF-T model to the basic MF model. We then show that the MF-MM model is a special case of the MF-T model.

\subsection{Mean Field Sharing with Multiple Types (MF-T) Model}\label{Multiple_types_model}
In this section, we describe the model of \textit{mean-field sharing with multiple types} (MF-T for short). This model is a generalization of the basic MF model. Instead of homogeneous subsystems, the subsystems are heterogeneous. In particular, subsystems are allowed to have different types such that, within each type, subsystems are homogeneous.

Let $\Theta$ denote the set of all different types of subsystems. Assume we have $k \in \mathbb{N}$ different types i.e. $\Theta=\{1,\ldots,k\}$. Let subsystem $i \in \{1,\ldots,n\}$ have a type  $\theta^i \in \Theta$. Let $X^i_t \in \mathcal{X}(\theta^i)$ denote the state of subsystem $i$ and $u^i_t \in \mathcal{U}(\theta^i)$ denote its control action. Assume $\mathcal{X}(\theta)$ and $\mathcal{U}(\theta)$ are finite sets, $\forall \theta \in \Theta$. Let $Z_t$ be as follows:
\begin{equation}\label{Multiple_types_mean_field}
 Z_t=\frac{1}{n} \sum_{i=1}^n \delta_{\theta^i, X^i_t}.
\end{equation}
For example, suppose there are $n=5$ subsystems with two different types i.e.  $\Theta=\{1,2\}$. The state space of type $\theta=1$ is $\mathcal{X}(1)=\{1,2\}$ and the state space of type $\theta=2$ is $\mathcal{X}(2)=\{2,3,4\}$. Subsystems $i=1,3$ are of type $1$ and subsystems $i=2,4,5$ are of type $2$. At time $t$, $x^1_t=2,x^2_t=3,x^3_t=2,x^4_t=4,x^5_t=2$, then $z_t(\theta=1,x=1)=0,z_t(\theta=1,x=2)=\frac{2}{5},z_t(\theta=2,x=2)=\frac{1}{5},z_t(\theta=2,x=3)=\frac{1}{5},z_t(\theta=2,x=4)=\frac{1}{5}$. Hence, mean-field  $z_t=(0,\frac{2}{5},\frac{1}{5},\frac{1}{5},\frac{1}{5})$.

The dynamics of subsystem $i$ is given by 
\begin{equation}\label{Multiple_types_dynamics}
X^i_{t+1}=f_t(\theta^i, Z_t,X^i_t,U^i_t,W^i_t),  i \in \{1,\ldots,n\},\theta^i \in \Theta
\end{equation}
where the primitive random variables $(X^1_1,\ldots,X^n_1,$ $\{W^1_t\}_{t=1}^T,\ldots,\{W^n_t\}_{t=1}^T)$ are mutually independent. The dynamics $f_t$ and the probability distribution of $X^i_1$ and $W^i_t$ may depend on  type $\theta^i \in \Theta$ but they do not directly depend on $i$. In other words, if subsystems $i$ and $j$ have the same type i.e. $\theta^i=\theta^j$, then their dynamics and probability distribution of their respective primitive random variables are the same. In addition, type of each subsystem is fixed and  every subsystem a priori knows the types of all other subsystems. 


Similar to the basic model, we assume that the control laws of all controllers with the same type are identical. Thus, the control actions are generated as follows:
\begin{equation}\label{Multiple_types_control_law}
\begin{aligned}
U^i_t=g_t(\theta^i, Z_{1:t},X^i_t)\quad i \in \{1,\ldots,n\},\theta^i \in \Theta.
\end{aligned}
\end{equation}

The collection $(\mathbf{g}(1), \ldots, \mathbf{g}(k))$, where $\mathbf{g}(\theta):=(g_1(\theta,\cdot),\ldots,g_T(\theta,\cdot),\theta \in \Theta$, is called the \textit{control strategy} of the system. Let $\mathcal{G}_{\mathit{MF-T}}$ be a set of all such control strategies.  At time $t$, the system incurs a cost depending on joint state $\mathbf X_t$ and joint action $\mathbf U_t$  that is given by
\begin{equation}\label{Multiple_types_cost}
\ell_t(\mathbf{X}_t,\mathbf{U}_t).
\end{equation}
The performance of any strategy  $(\mathbf{g}(1), \ldots, \mathbf{g}(k))$ is measured by the following expected total cost
\begin{equation}\label{Multiple_types_total_cost}
 J(\mathbf{g}(1),\ldots,\mathbf{g}(k))\hspace{-.1cm}=\hspace{-.1cm}\mathbb{E}^{\mathbf{g}(1),\ldots,\mathbf{g}(k)}  \left[\sum_{t=1}^T \hspace{-.1cm}\ell_t(\mathbf{X}_t,\mathbf{U}_t)\right]
\end{equation}
where the expectation is with respect to a joint measure induced on all system variables by the choice of $(\mathbf{g}(1), \ldots, \mathbf{g}(k))$. We are interested in the following optimization problem.

\begin{Problem}
For the MF-T model, identify a control strategy $({\mathbf{g}(1)}^\ast, \ldots, {\mathbf{g}(k)}^\ast)$ to minimize the expected total cost given by \eqref{Multiple_types_total_cost}.
\end{Problem}

\section{Main Results}\label{Main_results}

In this section, we present the main results of three mean-field models described in Section \ref{MF_models}.
\subsection{Basic MF Model}\label{BMF_main_results}

The basic MF model was investigated in \cite{Jalal2014MF}. The key idea in \cite{Jalal2014MF} is to construct a \textit{coordinated system} using the common information approach of \cite{Nayyar2013CIA} and show that the mean-field $Z_t$ is an information state for the coordinated system. Using this information state, we can obtain a dynamic program to identify an optimal control strategy. For completeness, we summarize the main steps below.

Following \cite{Nayyar2013CIA}, we construct a coordinated system from original system described in Section \ref{BMF_model}. Consider a virtual \textit{coordinator} that observes the common information between all subsystems by time $t$ i.e. $Z_{1:t}$. Based on $Z_{1:t}$,  the coordinator decides control action  $\Gamma_t: \mathcal{X} \rightarrow \mathcal{U}$ as follows:
\begin{equation}\label{Psi}
\Gamma_t =\psi_t(Z_{1:t}).
\end{equation}
We refer to $\psi_t$ as \textit{coordinator}'s control law at time $t$. In the coordinated system, each controller is a passive agent that uses its local state $X^i_t$ and the mapping $\Gamma_t$, which is called \textit{prescription}, to generate  
\begin{equation}\label{BMF_main_results_gamma}
U^i_t=\Gamma_t(X^i_t), \quad i \in \{1,\ldots,n\}.
\end{equation}
Note that $\Gamma_t$ and $\psi_t$ do not depend on $i$. The dynamics of each subsystem and the cost function are the same as in the original problem. 

It is shown in \cite{Nayyar2013CIA} that the coordinated system is equivalent to the original system. In particular, for any $\psi_t$, we can construct $g_t$ given by 
\begin{equation}
g_t(Z_{1:t},X^i_t):=\psi_t(Z_{1:t})(X^i_t)
\end{equation}
such that the expected total cost under strategy $\bm g$ in the original system is the same as that under strategy $\bm \psi:=(\psi_1,\ldots,\psi_T)$ in the coordinated system.

Thus, instead of searching for an optimal strategy $\bm g=(g_1,\ldots,g_T)$ in the original system, we may, equivalently, search for an optimal coordination strategy $\bm \psi=(\psi_1,\ldots,\psi_T)$ in the coordinated system.

It is shown in \cite{Jalal2014MF} that $Z_t$ is an information state for coordinated system. In particular, the cost given in \eqref{BMF_cost} may be viewed as a function of $(Z_t,\Gamma_t)$  (Lemma \ref{BMF_Z_cost}) and the evolution of mean-field $Z_{t+1}$ depends only on the current information $\{Z_t,\Gamma_t\}$ (Lemma \ref{BMF_Z_evolution}). 
 
\begin{Lemma}[Lemma 3 of \cite{Jalal2014MF}]\label{BMF_Z_cost}
The per-step cost may be written as a function of $Z_t$ and $\Gamma_t$. In particular,
there exits a function $\hat{\ell}$ such that 
$ \Exp{\ell_t(\mathbf{X}_t, \mathbf U_t)|Z_{1:t},\Gamma_{1:t}}=:\hat{\ell}_t(Z_t,\Gamma_t)$.
\end{Lemma}
\begin{Lemma}[Lemma 4 of \cite{Jalal2014MF}]\label{BMF_Z_evolution}
For any choice $\gamma_{1:t}$ of $\Gamma_{1:t}$, any realization $z_{1:t}$ of $Z_{1:t}$, and any $z$,
\[
\Prob{Z_{t+1} \hspace{-.1cm}= \hspace{-.1cm} z|Z_{1:t} \hspace{-.1cm} = \hspace{-.1cm} z_{1:t},\hspace{-.05cm}\Gamma_{1:t} \hspace{-.1cm} =\hspace{-.1cm} \gamma_{1:t}} \hspace{-.1cm}= \hspace{-.1cm} \Prob{Z_{t+1} \hspace{-.1cm} =\hspace{-.1cm}z|Z_t \hspace{-.1cm} = \hspace{-.1cm} z_{t}, \hspace{-.05cm}\Gamma_t \hspace{-.1cm} =\hspace{-.1cm} \gamma_{t}}.\]
\end{Lemma}
From Lemmas \ref{BMF_Z_cost} and \ref{BMF_Z_evolution}, one can conclude that $Z_t$ is  a controlled Markov process with control action $\Gamma_t$. Thus, standard results for MDPs can be used to obtain the optimal strategy as follows.

\begin{Theorem}[Theorem 1 of \cite{Jalal2014MF}]\label{BMF_theorem_complete}
Within basic mean-field control strategies $\mathcal{G}_{\mathit{BMF}}$, there is no loss of optimality in restricting attention to Markovian strategy i.e. $U^i_t=g_t(Z_t,X^i_t)$. Furthermore,  an optimal strategy $\bm g^\ast$  is obtained by  solving the following dynamic program. Define recursively value functions 
$V_{T+1}(z):= 0, \forall z,$
and for $t= T,\ldots,1,$
\begin{equation}\label{BMF_dp}
V_t(z_t)\hspace{-.1cm}:=\hspace{-.1cm}\min_{\gamma_t}(\hat{\ell}_t(z_t,\hspace{-.05cm} \gamma_t)\hspace{-.05cm}+ \hspace{-.05cm}\Exp{V_{t+1}(Z_{t+1})|Z_t\hspace{-.1cm}=\hspace{-.1cm}z_t,\hspace{-.05cm} \Gamma_t\hspace{-.1cm}=\hspace{-.1cm}\gamma_t})
\end{equation}
where $\gamma_t:\mathcal{X} \mapsto \mathcal{U}$. Let $ \psi_t^\ast(z_t)$ denote the argmin in the right-hand side of \eqref{BMF_dp}. Define $g^\ast_t(z,x):=\psi^\ast_t(z)(x), \forall x,z$. Then, $\mathbf g^\ast=(g^\ast_1,\ldots,g^\ast_T)$ is an optimal strategy.
\end{Theorem} 
In \cite{Jalal2014MF}, it is shown that the above results extend to the case when all subsystems observe a noisy version of the mean field by defining the information state as the posterior probability of the mean-field given noisy observations. 
\hspace{-.05cm}

\subsection{MF-T Model}\label{Multiple_types_main_results}

The MF-T model is a straightforward generalization of the basic MF model and can be solved in a similar manner as Section \ref{BMF_main_results}. For succinctness, we simply present a change of variable argument that lets us view the MF-T model as a special case of the basic MF model.

Consider the MF-T model in Section \ref{Multiple_types_model}. Define $\tilde{X}^i_t:=(\theta^i, X^i_t)$ as an augmented state that takes value in state space $\mathcal{\tilde{X}}= \cup_{\theta \in \Theta} \left( \{\theta\} \times \mathcal{X}(\theta)\right)$.  With this augmented state space, the MF-T model is converted to the basic MF model. In particular, mean-field of the MF-T model \eqref{Multiple_types_mean_field} may be written as 
$Z_t=\frac{1}{n} \sum_{i=1}^n \delta_{\tilde{X}^i_t}$
which is in the form of \eqref{BF_mean_field}. Similarly, the dynamics given by \eqref{Multiple_types_dynamics} may be written as
\begin{equation}
\tilde{X}^i_{t+1}=\tilde{f}_t( Z_t,\tilde{X}^i_t,U^i_t, W^i_t)
\end{equation}
which is of the form of \eqref{BMF_dynamics}. Since the type of subsystems are fixed and a priori known to every subsystem, the primitive random variables $(\tilde X^1_1,\ldots,\tilde X^n_1, \{W^1_t\}_{t=1}^T,\ldots,\{W^n_t\}_{t=1}^T)$ are mutually  independent conditioned on the types of subsystems. Note that  dynamics $\tilde{f}_t$ and probability distribution on $\tilde{X}^i_t$ and $W^i_t$ may only depend on the types and do not directly depend on the index $i$.

In addition, control law of the MF-T model \eqref{Multiple_types_control_law} may be formulated as
\begin{equation}
 U^i_t=\tilde{g}_t(Z_{1:t},\tilde{X}^i_t)
\end{equation}
which is in the form of \eqref{BMF_control_law}. Furthermore, the cost given by \eqref{Multiple_types_cost} may be written as $\tilde{\ell}_t(\mathbf{\tilde{X}_t}, \mathbf U_t)$
which is of the form of \eqref{BMF_cost}. 
Thus, similar to Section \ref{BMF_main_results}, we define function $\tilde{\psi}_t$ (coordinator's control law) as a mapping from the history of mean-field  $ Z_{1:t}$ to $\tilde \Gamma_t$ in the same fashion as \eqref{Psi} such that 
$\tilde{\Gamma}_t=\tilde \psi_t(Z_{1:t})$
where 
$U^i_t=\tilde{\Gamma}_t(\tilde X^i_t),  i \in \{1,\ldots,n\},$
which is of the form of \eqref{BMF_main_results_gamma}.

Hence, an immediate consequence of the results of Section \ref{BMF_main_results} are the following:

\begin{Lemma}\label{Multiple_types_Z_cost}
 The per-step cost given in $\tilde{\ell}_t(\mathbf{\tilde{X}_t}, \mathbf U_t)$ may be written as a function of $Z_t$ and $\tilde \Gamma_t$. In particular,
there exits a function $\hat{\ell}$ such that 
$ \Exp{\tilde{\ell}_t(\mathbf{\tilde{X}}_t, \mathbf U_t)|Z_{1:t},\tilde{\Gamma}_{1:t}}=:\hat{\ell}_t(Z_t,\tilde{\Gamma}_t)$.
\end{Lemma}
\begin{Lemma}\label{Multiple_types_Z_evolution}
For any choice $\tilde \gamma_{1:t}$ of $\tilde \Gamma_{1:t}$, any realization $ z_{1:t}$ of $Z_{1:t}$, and any $ z$,
\[
\mathbb{P}(\hspace{-.05cm} Z_{t+1} \hspace{-.1cm}= \hspace{-.1cm}  z| Z_{1:t} \hspace{-.1cm} = \hspace{-.1cm}  z_{1:t},\tilde \Gamma_{1:t} \hspace{-.1cm} =\hspace{-.1cm} \tilde \gamma_{1:t}\hspace{-.05cm} ) \hspace{-.1cm}= \hspace{-.1cm}  \Prob{ Z_{t+1} \hspace{-.1cm} = \hspace{-.1cm}z| Z_t \hspace{-.1cm} = \hspace{-.1cm}  z_{t}, \tilde \Gamma_t \hspace{-.1cm} =\hspace{-.1cm} \tilde \gamma_{t}}.\]
\end{Lemma}
Similar to Section \ref{BMF_main_results}, based on Lemmas \ref{Multiple_types_Z_cost} and \ref{Multiple_types_Z_evolution}, one can conclude that $Z_t$ is  a controlled Markov process with control action $\tilde \Gamma_t$. Thus, standard results for MDPs can be used to obtain the optimal strategy as follows.
\begin{Theorem}\label{Multple_types_theorem_complete}
Within MF-T control strategies $\mathcal{G}_{\mathit{MF-T}}$, there is no loss of optimality in restricting attention to Markovian strategy i.e. $U^i_t=g_t(\theta^i, Z_t, X^i_t)$. Furthermore, an optimal strategy $(\bm g(1)^\ast,\ldots, \bm g(k)^\ast)$ is obtained by solving the following dynamic program. Define recursively value functions 
$V_{T+1}(z):= 0, \forall z,$
and for $t= T,\ldots,1,$
\begin{equation}\label{Multiple_types_dp}
V_t(z_t)\hspace{-.1cm}:=\hspace{-.1cm}\min_{\tilde \gamma_t}(\hat{\ell}_t(z_t,\hspace{-.05cm}\tilde \gamma_t)\hspace{-.04cm}+ \hspace{-.04cm}\Exp{V_{t+1}(\hspace{-.05cm}Z_{t+1}\hspace{-.05cm})|Z_t\hspace{-.1cm}=\hspace{-.1cm}z_t,\hspace{-.05cm} \tilde \Gamma_t\hspace{-.1cm}=\hspace{-.1cm}\tilde \gamma_t})
\end{equation}
where $\tilde \gamma_t=(\gamma_t(1),\ldots,\gamma_t(k))$ and $\gamma_t(\theta):\mathcal{X}(\theta) \mapsto \mathcal{U}(\theta), \forall \theta \in \Theta$.  Let $ \psi_t^\ast(z_t)$ denote the argmin in the right-hand side of \eqref{Multiple_types_dp}. Define $g^\ast_t(\theta,z,x):=\psi^\ast_t(z)(\theta,x), \forall \theta,x,z$. Then,   the collection of strategies ${\mathbf g(\theta)}^\ast=({g_1(\theta,\cdot)}^\ast,\ldots,{g_T(\theta,\cdot)}^\ast), \forall \theta \in \Theta$, is  optimal.
\end{Theorem}

\subsection{MF-MM Model}\label{MMMF_main_results}

In this section, we show that the MF-MM model is a special case of the MF-T model described in Section \ref{Multiple_types_model}. In particular, the system has two types: major and minor. There is only one major subsystem; hence the mean-field is equivalent to $(Z_t,X^0_t)$. Let $\Gamma_t: \mathcal{X} \mapsto \mathcal{U}$ and $\Gamma^0_t: \mathcal{X}^0 \mapsto \mathcal{U}^0$. Then, by virtue of Lemma \ref{Multiple_types_Z_cost} and Lemma \ref{Multiple_types_Z_evolution}, we have that $(Z_t,X^0_t)$ is a controlled Markov process with control action $(\Gamma_t,\Gamma^0_t)$. Note that this result is  simplified further because $X^0_t$ is part of the local information at the major subsystem as well as the common information. Therefore, it can be shown that, the function $\Gamma^0_t: \mathcal{X}^0 \mapsto \mathcal{U}^0$ may be replaced by the action $U^0_t$. Thus, from Lemma \ref{Multiple_types_Z_cost} and Lemma \ref{Multiple_types_Z_evolution}, we get the following:

\begin{Lemma}\label{MMMF_Z_cost}
The per-step cost $ \ell_t(\mathbf{X}_t, X^0_t, \mathbf{U}_t, U^0_t)$ may be written as a function of $(Z_t,X^0_t)$ and $(\Gamma_t,U^0_t)$. In particular,
there exits a function $\hat{\ell}$ such that 
\[ \Exp{\ell_t\hspace{-0.05cm}(\mathbf X_t,\hspace{-0.05cm}X^0_t,\hspace{-0.05cm}\mathbf U_t,\hspace{-0.05cm}U^0_t)|Z_{1:t},\hspace{-0.05cm}X^0_{1:t},\hspace{-0.05cm}\Gamma_{1:t},\hspace{-0.05cm}U^0_{1:t}}\hspace{-0.1cm}=:\hspace{-0.1cm}\hat{\ell}_t\hspace{-0.05cm}(Z_t,\hspace{-0.05cm}X^0_t,\hspace{-0.05cm}\Gamma_t,\hspace{-0.05cm}U^0_t).\]
\end{Lemma}
\begin{Lemma}\label{MMMF_Z_evolution}
For any choice of $(\gamma_{1:t}, u^0_{1:t})$ of
 $(\Gamma_{1:t}, U^0_{1:t})$, any realization $(z_{1:t},x^0_{1:t})$ of $(Z_{1:t},X^0_{1:t})$, and any $(z,x^0)$,
\begin{equation*}
\begin{aligned}
&\hspace*{-.1cm}\mathbb{P}\hspace{-.03cm}(\hspace{-.02cm} Z_{t+1}\hspace{-.15cm}=\hspace{-.1cm}z, \hspace{-.05cm}X^0_{t+1}\hspace{-.15cm}=\hspace{-.1cm}x^0 \hspace{-0.02cm}|\hspace{-0.02cm} Z_{1:t}\hspace{-.1cm}=\hspace{-.1cm}z_{1:t} ,\hspace{-.05cm}X^0_{1:t}\hspace{-.1cm}=\hspace{-.1cm}x^0_{1:t},\hspace{-.05cm} \Gamma_{1:t}\hspace{-.1cm}=\hspace{-.1cm}\gamma_{1:t},\hspace{-.05cm} U^0_{1:t}\hspace{-.1cm}=\hspace{-.1cm}u^0_{1:t}\hspace{-.05cm})\\
&\hspace{-.2cm}=\mathbb{P}\hspace{-.03cm}(\hspace{-.03cm} Z_{t+1}\hspace{-.1cm}=\hspace{-.1cm}z ,X^0_{t+1}\hspace{-.1cm}=\hspace{-.1cm}x^0 | Z_{t}\hspace{-.1cm}=\hspace{-.1cm}z_{t} ,X^0_{t}\hspace{-.1cm}=\hspace{-.1cm}x^0_{t}, \Gamma_t\hspace{-.1cm}=\hspace{-.1cm}\gamma_t, U^0_t\hspace{-.1cm}=\hspace{-.1cm}u^0_t)
\end{aligned}
\end{equation*}
\end{Lemma}
From  Lemmas \ref{MMMF_Z_cost} and \ref{MMMF_Z_evolution}, we can conclude that $(Z_t,X^0_t)$ is a controlled Markov process with control actions $(\Gamma_t,U^0_t)$. Hence,  we can use the standard results for MDPs to obtain the optimal strategies.
\begin{Theorem}
Within MF-MM control strategies $\mathcal{G}_{\mathit{MF-MM}}$, there is no loss of optimality in restricting attention to Markovian strategy i.e. $U^0_t=g^0_t(Z_t,X^0_t)$ and $U^i_t=g_t(Z_t,X^0_t,X^i_t)$. Furthermore, an optimal strategy $({\bm {g}^0}^\ast,\bm g^\ast)$ is obtained by solving the following dynamic program. Define recursively value functions 
$V_{T+1}(z,x^0):= 0, \forall z,x^0,$
and for $t= T,\ldots,1,$
\begin{equation}\label{MF-MM_dp}
\begin{aligned}
&V_t(z_t,x^0_t):=\min_{\gamma_t,u^0_t}(\hat{\ell}_t(z_t,x^0_t,\gamma_t,u^0_t)\\
& \hspace*{-.3cm}+ \hspace{-.05cm}\Exp{V_{t+1}\hspace{-.05cm}(\hspace{-.05cm}Z_{t+1},\hspace{-.05cm}X^0_{t+1}\hspace{-.05cm})\hspace{-.02cm}|\hspace{-.03cm}Z_t\hspace{-.1cm}=\hspace{-.1cm}z_t,\hspace{-.05cm}X^0_t\hspace{-.1cm}=\hspace{-.1cm}x^0_t,\hspace{-.05cm}\Gamma_t\hspace{-.1cm}=\hspace{-.1cm}\gamma_t,\hspace{-.05cm} U^0_t=u^0_t})
\end{aligned}
\end{equation}
where $\gamma_t:\mathcal{X} \mapsto \mathcal{U}$. Let $ \psi_t^\ast(z_t,x^0_t)$ denote the argmin in the right-hand side of \eqref{MF-MM_dp} such that $\gamma_t^\ast=  \psi_t^\ast(1)(z_t,x^0_t)$ and ${u_t^0}^\ast=\psi_t^\ast(2)(z_t,x^0_t)$. Define $g^\ast_t(z,x^0,x):={\psi^\ast_t(1)}(z,x^0)(x)$ and ${g^0_t}^\ast(z,x^0):={\psi^\ast_t(2)}(z,x^0), \forall x,z,x^0$. Then, the collection of strategies $\mathbf g^\ast=(g^\ast_1,\ldots,g^\ast_T)$ and ${\mathbf{g}^0}^\ast=({g^0_1}^\ast,\ldots,{g^0_T}^\ast)$ is  optimal.
\end{Theorem}
\begin{Remark}
At time $t$, the optimal action of major subsystem ${U^0_t}^\ast$ depends on the mean field of minor subsystems $Z_t$ and the state of major subsystem $X^0_t$ i.e. ${U^0_t}^\ast={g^0_t}^\ast(Z_t,X^0_t)$. Moreover, the optimal action of minor subsystem $i$, ${U^i_t}^\ast$, depends on its own local state $X^i_t$ in addition to   $Z_t$ and $X^0_t$ i.e. ${U^i_t}^\ast=g^\ast_t(Z_t,X^0_t,X^i_t)$. 
\end{Remark}

\section{Example}\label{Example}
In this section, we present an example of the MF-MM model where a service provider (major subsystem) interacts with its users (minor subsystems). Roughly speaking, we answer the question of how a service provider should design strategies for itself and the users such that the service is not only profitable but also customer-satisfactory.

Consider a service provider (internet provider, utility provider, etc.) who serves $ n \in \mathbb{N}$ users (customers). 
Let $X^i_t \in \mathcal{X}=\{0,1\}$, $i \in \{1,\ldots,n\}$, denote whether  user $i$ is active $(X^i_t=1)$ or passive $(X^i_t=0)$ at time $t$ and $X^0_t \in \mathcal{X}^0$ denote the available capacity of the system. At time $t$, the service provider chooses $U^0_t$, the capacity of the system at time $t+1$. Thus, the dynamics of the state of the service provider is
\begin{equation}
X^0_{t+1}=U^0_t
\end{equation}
where $\mathcal{X}^0=\mathcal{U}^0$ is a finite set consisting of all feasible capacities.

Dynamics of users are represented by $P(u)$ where $P(u)$ denotes the controlled transition matrix under action $u \in \mathcal{U}=\{0,1,2\}$, i.e.
\begin{equation*}
[P(u)]_{xy}=\Prob{X^i_{t+1}=y \mid X^i_t=x, U^i_t=u},\quad x,y \in \mathcal{X}.
\end{equation*} 
The service provider can affect the evolution of user $i$ by choosing action $u^i \in \{0,1,2\}$. Action $u^i=0$ is a \textit{free action} under which user $i$ evolves in an uncontrolled manner, i.e. $P(0)=Q$, where $Q$ represents the \textit{natural} dynamics of the users. Action $u^i \neq 0$ is a \textit{forcing action} under which service provider forces user $i$ to switch to state $(u^i-1)$. However,  there is $\epsilon_{u^i} \in [0,1]$ probability that user $i$ does \emph{not} switch to state $(u^i-1)$ and  evolves according to its natural dynamics. Thus,
\begin{equation*}
P(u)=(1-\epsilon_u)\mathbf K_u + \epsilon_u Q
\end{equation*}
where $\mathbf K_u$ is a $k \times k$ matrix where column $u \in \{1,2\}$ is all ones, and other columns are all zeros.
%
%
%
The goal is that the service provider efficiently manages the users and the capacity of the system in order to minimize costs while satisfying the users. Thus, the loss function is given by
\begin{equation}
\ell_t^0(Z_t,X^0_t,U^0_t)+ \sum_{i=1}^n \ell_t(U^i_t)
\end{equation}
where $\ell_t^0$ is as follows:
\begin{equation}
\ell_t^0(Z_t,X^0_t,U^0_t)=S(U_t^0) + a|U^0_t-X^0_t| - G(Z_t,X^0_t), \hspace{.1cm} a \in \mathbb{R}^+
\end{equation}
and 
$\ell_t(U^i_t)=H(U^i_t).$
In $\ell^0_t$, the first term is associated with the cost of capacity and the second term refers to patching and dispatching capacity. The third term corresponds to the benefit (proportional to the number of active users) and the penalty of unavailable service (proportional to the number of active users that do not receive service) i.e.
\begin{equation}
G(Z_t,X^0_t)=
\begin{dcases}
bnZ_t(1) & nZ_t(1)\leq X^0_t\\
bX^0_t-c(nZ_t(1)-X^0_t)  & nZ_t(1)> X^0_t
\end{dcases}
\end{equation}
where $b \in \mathbb{R}^+$ indicates the rate of benefit and $c \in \mathbb{R}^+$ determines the rate of penalty of unavailable service. Note that $Z_t(1)=\frac{1}{n}\sum_{i=1}^n \mathds{1}(X^i_t=1)$ is the average number of active users  at time $t$. In addition, $\ell_t$ is the cost associated with forcing users.

Given the information structure \eqref{MF-MM_information_structure}, the objective is to choose a control strategy that minimizes the infinite horizon discounted cost
\begin{equation}\label{Example_J}
J(\bm g)\hspace{-.1cm}= \mathbb{E}\hspace{-.1cm}\left[\sum_{t=1}^\infty \hspace{-.05cm}\beta^t \hspace{-.1cm}\left( \ell_t^0(Z_t,X^0_t,U^0_t)+\hspace{-.1cm} \sum_{i=1}^n \hspace{-.1cm} \ell_t(U^i_t)\right)\hspace{-.05cm}\right]
\end{equation}
where $\beta \in (0,1)$ is a discount factor\footnote{Although we have only presented the details for finite horizon setup in this paper, the results generalize naturally to infinite horizon setup under standard assumption.}.
The optimal time-homogeneous strategies for 
\begin{equation*}
\begin{aligned}
&n=100, \quad \hspace{-.15cm} \mathcal{X}^0=\{50,100\}, \quad\hspace{-.15cm} S(50)=100,\quad \hspace{-.15cm} S(100)=300 ,\\
&H(0)=0, \quad \hspace{-.15cm} H(1)=4, \hspace{-.15cm} \quad H(2)=1, \quad \hspace{-.15cm}\epsilon_1 =0.1, \quad \hspace{-.15cm} \epsilon_2 =0.1,\\
&a=2,\quad b=5,\quad c=50, \quad \beta=0.6,  \quad 
  \quad \hspace{-.2cm} Q=\left[\begin{array}{cc} 
 0.6 \quad  0.4 \\  0.3 \quad  0.7 \end{array} \right] 
\end{aligned}
\end{equation*}
are presented below. Since the state space of minor subsystems is binary, $z(1)$ is sufficient to characterise the empirical distribution $z=[z(0), z(1)]$. Hence, for ease of presentation, we represent the optimal control law as a function of the second component $z(1)$ of $z=[z(0), z(1)]$.
\begin{equation*}\label{Minor_solution}
g^\ast(z,x^0,x)=
\begin{dcases}
\begin{aligned}
0 & &0\leq z(1)\leq 0.53&,&  x^0=50&,  x=0\\
1 & &0.53< z(1)\leq 0.76&,& x^0=50&, x=0\\
2 & &0.76< z(1)\leq 1&,&  x^0=50&, x=0\\
0 & &0\leq z(1)\leq 1&,& x^0=50&, x=1\\
0 & &0\leq z(1)\leq 0.29&,& x^0=100&, x=0\\
2 & &0.29< z(1)\leq 1&,&  x^0=100&, x=0\\
0 & &0\leq z(1)\leq 1&,&  x^0=100&, x=1\\
\end{aligned}
\end{dcases}
\end{equation*}
and
\begin{equation*}\label{Major_solution}
{g^0}^\ast(z,x^0)=
\begin{dcases}
\begin{aligned}
50 & &0\leq z(1)\leq 0.76,& \quad x^0=50\\
100 & &0.76< z(1)\leq 1,& \quad  x^0=50\\
50 & &0\leq z(1)\leq 0.29,&\quad  x^0=100\\
100 & &0.29< z(1)\leq 1,& \quad  x^0=100.\\
\end{aligned}
\end{dcases}
\end{equation*}
%
\section{Conclusion}\label{Conclusion}
 In this paper, we presented  two variations of the basic mean-field sharing model: mean-field sharing  with major-minor subsystems (MF-MM) and  mean-field sharing  with multiple types (MF-T). We demonstrated that, by an appropriate change of variables, the MF-T model may be viewed as a special case of the basic mean-field (MF) model presented in~\cite{Jalal2014MF}. Based on the main results of the basic MF model in~\cite{Jalal2014MF},  we identified a dynamic program to identify optimal strategies for the MF-T model. In addition, we showed that the MF-MM model may be viewed as a special case of the MF-T model. Due to the special structure of the MF-MM model, we could simplify the dynamic program associated with the MF-T model further  to obtain the optimal strategies for the major and  minor subsystems. Note that we provided a straightforward solution  for the MF-MM model in contrast to mean-field games where such a solution is non-trivial and difficult. 

\vspace{-.0cm}
\bibliography{ICC2015_Ref}
\bibliographystyle{IEEEtran}
\end{document}